\documentclass[12pt,twoside,a4paper]{article}
\usepackage{srcltx}
\usepackage{amsmath,amsfonts,amscd,amssymb,latexsym}
\usepackage{xcolor}
\usepackage[latin1]{inputenc}
\newtheorem{theorem}{Theorem}[section]
\newtheorem{lemma}{Lemma}[section]
\newtheorem{corollary}{Corollary}[section]
\newtheorem{remark}{Remark}[section]
\newcommand{\oti}[2]{\mathop \otimes\limits_{#2}^{#1} v}

\newcommand{\ox}[3]{\mathop {#1}\limits_{#3}^{#2} }

\newcommand{\0}{_{\bar 0}}

\newcommand{\1}{_{\bar 1}}
\newcommand{\F}{{\cal F}}

\def\proof{{\noindent\it Proof. }}
\def\qed{\hfill \rule{2.25mm}{2.25mm}\vspace{0,2cm}}
\renewcommand{\=}{\doteq}
\def\qed{\hfill \rule{2.25mm}{2.25mm}\vspace{10pt}}
\renewcommand{\L}{{\cal L}}

\begin{document}

\noindent
\vspace{1 in}
\begin{center}
\noindent {\Large {\bf ON SIMPLE FILIPPOV
SUPERALGEBRAS\\[2mm] OF TYPE $A(m,n)$}}
\vspace{.5cm}

P. D. Beites and A. P. Pozhidaev
\end{center}
\vspace{0,5cm}

\parbox[c]{13cm}{{\large Abstract}$:$
It is proved that there exist no simple finite-dimensional
Filippov superalgebras of type $A(m,n)$ over an
algebraically closed field of characteristic $0$.\\

}

\section{Introduction}

  The concept of $n$-Lie superalgebra was
presented by Daletskii and Kushnirevich, in \cite{Dal}, as a natural generalization of the
  $n$-Lie algebra notion introduced by Filippov in $1985$ (see \cite{Fil}).
  Following \cite{GM} and \cite{EnL}, in this article, we use the terms
  Filippov superalgebra and Filippov algebra
  instead of $n$-Lie superalgebra and  $n$-Lie algebra, respectively. Filippov
  algebras were also known before under the names of
 Nambu Lie algebras and Nambu algebras. As pointed out in \cite{malc} and \cite{genmalc}, Filippov algebras are
 a particular case of $n$-ary Malcev algebras (generalizing the fact that every Lie algebra is a Malcev algebra).
 We may also remark that a $2$-Lie superalgebra is simply known as a Lie superalgebra. The description of the
 finite-dimensional simple Lie superalgebras over an algebraically closed field of characteristic zero was given
 by Kac in \cite{Kac}.

This work is one more step on the way to the classification of
  finite-dimensional simple Fi\-lippov superalgebras over an algebraically
  closed field of characteristic $0$.
  In \cite{LA}, finite-dimensional commutative $n$-ary Leibniz
  algebras over a field of characteristic $0$ were studied by the second author. He showed
  that there exist no simple ones.
  The finite-dimensional simple Filippov algebras over an
  algebraically closed field of characteristic $0$ were classified by Ling in
  \cite{Ling}. Notice that an $n$-ary commutative Leibniz algebra is
  exactly a Filippov superalgebra with trivial even part, and a Filippov
  algebra is exactly a Filippov superalgebra with trivial odd part.
Bearing in mind these facts, we consider the
 $n$-ary Filippov superalgebras
  with $n\geq 3$, and with nonzero even and odd parts.

Let $G$ be a Lie superalgebra. We say that a Filippov superalgebra
$\F$ has {\it type} $G$ if $Inder(\F)\cong G$ (see definitions below).
A description of simple Filippov superalgebras of type $B(m,n)$ was already
obtained in \cite{Poj8}, \cite{PS} and \cite{P7}.
 The same problem concerning Filippov superalgebras
of type $A(m,n)$ with $m=n$ has recently been solved in \cite{BPoj}.
Moreover, the type $A(0,n)$, with $n \in \mathbb{N}$, was studied in \cite{pat_artigo4}.
The present work represents the final step towards
the classification of finite-dimensional simple Fi\-lippov superalgebras of type
$A(m,n)$ over an algebraically closed field of characteristic zero. Concretely, we
establish a negative answer to the existence problem of the mentioned superalgebras
 when $m, n \in \mathbb{N}$ and $m \neq n$.

We start recalling some definitions.

 An {\it $\Omega$-algebra} over a field $k$ is a linear
  space over $k$ equipped with a system of multilinear
  algebraic operations $\Omega = \{\omega_i: | \omega_i | = n_i \in
  {\mathbb N}, \ i \in I \}$, where $| \omega_i |$
  denotes the arity of $\omega_i$.

  An {\it $n$-ary Leibniz algebra} over a field $k$ is
  an $\Omega$-algebra $L$ over $k$ with one
  $n$-ary operation $(\cdot, \cdots, \cdot)$ satisfying the identity
 \begin{center}
  $((x_1, \ldots, x_n), y_2, \ldots, y_n) =
  \displaystyle\sum_{i = 1}^{n} (x_1, \ldots, (x_i, y_2, \ldots, y_n), \ldots, x_n).$ \end{center}  If this operation is anticommutative, we obtain the definition of
  \textit{Filippov} (\textit{$n$-Lie}) algebra over a field.

  An  {\it $n$-ary superalgebra} over a field $k$ is
  a $\mathbb{Z}_2$-graded $n$-ary algebra $L=L\0\oplus L\1$ over $k$,
  that is,

\begin{center} if $x_i\in L_{\alpha_i}, \alpha_i\in \mathbb{Z}_2$, then
  $(x_1, \ldots, x_n)\in L_{\alpha_1+\ldots+\alpha_n}$. \end{center}
An {\it $n$-ary Filippov superalgebra} over $k$ is
  an $n$-ary superalgebra $\F=\F\0\oplus\F\1$ over $k$ with one
  $n$-ary operation $[\cdot, \cdots, \cdot]$ satisfying
 \begin{eqnarray}\label{AC}
  &&\!\!\!\! \!\!\!\![x_1, \ldots, x_{i-1},x_i, \ldots, x_n] =
  -(-1)^{p(x_{i-1})p(x_i)}[x_1, \ldots, x_{i},x_{i-1}, \ldots, x_n],\\
  \label{YaI}
&&\!\!\!\!\!\!\!\![[x_1, \ldots, x_n], y_2, \ldots, y_n]  =
  \sum_{i = 1}^{n}(-1)^{p\bar q_i}
   [x_1, \ldots, [x_i, y_2, \ldots, y_n], \ldots, x_n],
 \end{eqnarray}
 where $p(x)=l$ means that $x\in \F_{\bar l}$, $p=\sum_{i=2}^{n}p(y_i),\
  \bar q_i=\sum_{j=i+1}^{n}p(x_j),\ \bar{q}_n=0$. The identities (\ref{AC})
  and (\ref{YaI})
 are called
 the anticommutativity and the generalized Jacobi identity, respectively.
 By (\ref{AC}), we can rewrite (\ref{YaI}) as
\begin{eqnarray}\label{LYaI}
&&\!\!\!\!\!\!\!\![y_2, \ldots, y_n,[x_1, \ldots, x_n]] =
  \sum_{i = 1}^{n}(-1)^{p q_i}
   [x_1, \ldots, [y_2,\ldots, y_n,x_i], \ldots, x_n],
 \end{eqnarray}
where $q_i=\sum_{j=1}^{i-1}p(x_j),\ q_1=0$. Sometimes, instead of using the long term
 ``$n$-ary superalgebra'', we simply say for short ``superalgebra''.
 If we denote by $L_x=L_{(x_1,\ldots,x_{n-1})}$
 the operator of left multiplication
 $L_xy=[x_1,\ldots,x_{n-1},y]$, then, by (\ref{LYaI}), we get
\begin{eqnarray*}\label{OLYaI}
&&[L_y,L_x] =
  \sum_{i = 1}^{n-1}(-1)^{p q_i}
   L(x_1, \ldots,L_yx_i, \ldots, x_{n-1}),
 \end{eqnarray*}
where $L_y$ is an operator of left multiplication and $p$ its parity.
(Here and afterwards, we denote the supercommutator by $[\,,]$).

Let $L=L\0\oplus L\1$ be an $n$-ary anticommutative superalgebra.
A \textit{subsuperalgebra} $B=B\0\oplus B\1$ of the superalgebra $L$,
$B_{\bar i}\subseteq L_{\bar i}$,
is a $\mathbb{Z}_2$-graded vector subspace of $L$ such that $[B, \ldots, B] \subseteq B$.
A subsuperalgebra $I$ of $L$ is called an \it ideal\/ \rm if
$[I,L, \ldots ,L] \subseteq I$.
The subalgebra (in fact, an ideal) $L^{(1)}=[L, \ldots ,L]$ of $L$
is called the \it derived subsuperalgebra\/ \rm of $L$.
Put $L^{(i)}=[L^{(i-1)}, \ldots ,L^{(i-1)}]$, $i\in {\mathbb N}, i>1$.
The superalgebra $L$ is called {\it solvable} if $L^{(k)}=0$ for some $k$.
Denote by $R(L)$ the maximal solvable ideal of $L$ (if it exists). If $R(L)=0$,
 the superalgebra $L$ is called {\it semisimple}.
The superalgebra $L$ is called \it simple\/ \rm if
$L^{(1)} \neq 0$ and $L$ lacks ideals other than $0$ or $L$.

 The article is organized as follows.

In the second section we recall how to reduce the classification problem
of simple Filippov superalgebras to some question about Lie superalgebras,
using the same ideas as in \cite{Ling}. Concretely, we consider an existence problem for some skewsymmetric homomorphisms of semisimple
Lie superalgebras and their faithful irreducible modules. This section is
followed with the third one where we collect some definitions and results on Lie superalgebras that we will apply in the two last sections. We also fix some notations with the same purpose.

The fourth section is devoted to the problem of existence of finite-dimensional
simple Filippov superalgebras of
type $A(m,n)$ with $m \neq n$.
We start with the particular case $A(1,n)$ in the first subsection, where, taking into account \cite{BPoj}
 and \cite{pat_artigo4}, it is assumed that
 $n \in \mathbb{N}\setminus \{1\}$. The main result of this article (Theorem \ref{mainresult}) is stated
  and proved in the second subsection.

In each of the two mentioned subsections we restrict our considerations to the case of the Lie
superalgebra that gives the name to the type and solve the existence problem of the mentioned skewsymmetric
homomorphisms. It turns out that the required homomorphisms do not exist.
Therefore, there are no simple finite-dimensional Filippov superalgebras
of type  $A(m,n)$ over an algebraically closed field of characteristic $0$.
Moreover, as a corollary of its proof, we see
that there is no simple finite-dimensional Filippov superalgebra ${\cal F}$ of type
$A(m,n)$ such that ${\cal F}$ is a highest weight module over $A(m,n)$.

In what follows, by $\Phi$ we denote  an algebraically closed field of
 characteristic $0$, by $F$ a field of characteristic $0$, by $k$
 a field and by $\left< w_{\upsilon}; \ \upsilon \in \Upsilon \right>$
  a linear space
spanned by the family of vectors
 $\{w_{\upsilon}; \ \upsilon \in \Upsilon \}$ over a field (the field is clear from the context).
The symbol $:=$ denotes an
 equality by definition.


\section{Reduction to Lie superalgebras}

From now on,  we denote by $\F$ an $n$-ary Filippov superalgebra.
Let us denote by $\F^*$ ($L(\F)$) the associative (Lie) superalgebra
generated  by the operators $L(x_1,\ldots,x_{n-1})$, $x_i\in \F$.
The algebra $L(\F)$ is called {\it the algebra of multiplications} of $\F$.

\begin{lemma} \textnormal{\cite{Poj8}} \it Let $\F=\F\0\oplus\F\1$ be a simple
finite-dimensional Filippov superalgebra over a field of characteristic $0$
with  $\F\1\neq 0$. Then $L=L(\F)=L\0\oplus L\1$ has nontrivial even and odd
parts.
\end{lemma}

\begin{theorem} \textnormal{\cite{Poj8}} \it If $\F$ is a simple
finite-dimensional Filippov superalgebra over a field of
characteristic $0$, then $L=L(\F)$ is a semisimple Lie superalgebra.
\label{teo1}
\end{theorem}

Given an $n$-ary superalgebra $A$ with a multiplication
$(\cdot,\cdots,\cdot)$, we have $End(A)=End\0A\oplus End\1A$. The element $D\in
End_{\bar s}A$ is called a {\it derivation} of degree $s$ of $A$ if, for every
$a_1,\ldots,a_n\in A, p(a_i)=p_i$, the following equality holds
$$D(a_1,\ldots,a_n)=\sum_{i=1}^{n}(-1)^{sq_i}(a_1,\ldots,Da_i,\ldots,a_n),$$
where $q_i=\sum_{j=1}^{i-1}p_j$.
We denote by $Der_{\bar s}A \subset End_{\bar s}A$
the subspace of all derivations
of degree $s$ and set $Der(A)=Der\0A\oplus Der\1A$.
The subspace $Der(A)\subset End(A)$ is easily seen to be closed under the
bracket

\begin{equation*} [a,b]=ab-(-1)^{deg(a)deg(b)}ba \end{equation*} (known as the \textit{supercommutator}) and it is called
{\it the superalgebra of derivations} of $A$.

\sloppy
Fix $n-1$ elements $x_1,\ldots,x_{n-1}\in A$, $i\in \{1, \ldots, n\}$,
and define a transfor\-mation
$ad_i(x_1,\ldots,x_{n-1})\in End(A)$ by the rule

\begin{equation*} ad_i(x_1,\ldots,x_{n-1})\, x=(-1)^{pq_i}(x_1,\ldots,x_{i-1},x,
x_{i},\ldots,x_{n-1}), \end{equation*} where $p=p(x),p_i=p(x_i),q_i=\sum_{j=i}^{n-1}p_j$.

\sloppy
If, for all $i=1,\ldots,n$ and $x_1,\ldots,x_{n-1}\in A$, the transformations
$ad_i(x_1,\ldots,x_{n-1})\in End(A)$ are derivations of $A$,
then we call them {\it strictly inner derivations} and $A$
 {\it an inner-derivation superalgebra $({\cal ID}$-superalgebra}).
Notice that the $n$-ary Fi\-lippov superalgebras and the $n$-ary commutative
Leibniz algebras are examples of ${\cal ID}$-superalgebras.

Now let us denote by $Inder(A)$ the linear space
 spanned by the strictly inner derivations of $A$.
 If $A$ is an $n$-ary $\cal ID$-superalgebra then it is easy to see
 that $Inder(A)$ is an ideal of $Der(A)$.

\begin{lemma} \textnormal{\cite{Poj8}} Given  a simple $\cal ID$-superalgebra $A$ over $k$,
the Lie superalgebra $Inder(A)$ acts faithfully and irreducibly on $A$.
\end{lemma}

Let $\F$ be an $n$-ary Filippov superalgebra over $k$.
Notice that the map $ad:=ad_n:\otimes^{n-1}\F\mapsto Inder(\F)$ satisfies
$$[D,ad(x_1,\ldots,x_{n-1})]=\sum_{i=1}^{n-1}(-1)^{pq_i}
ad(x_1,\ldots,x_{i-1},Dx_i,x_{i+1},\ldots,x_{n-1}),$$
for all $D\in Inder(\F)$, and the associated map
$(x_1,\ldots,x_n)\mapsto ad(x_1,\ldots,x_{n-1})\, x_n$
from $\otimes^{n}\F$ to $\F$ is $\mathbb{Z}_2$-skewsymmetric.
If we regard $\F$ as an $Inder(\F)$-module
then $ad$ induces an $Inder(\F)$-module morphism from
the $(n-1)$-th exterior power $\wedge^{n-1}\F$ to
$Inder(\F)$ (which we also denote by $ad$) such that the map
 $(x_1,\ldots,x_n)\mapsto ad(x_1,\ldots,x_{n-1})\, x_n$ is $\mathbb{Z}_2$-skewsymmetric.
(Note that in $\wedge^{n-1}\F$ we
have $x_1\wedge\ldots\wedge x_i\wedge x_{i+1}\wedge\ldots \wedge x_{n-1}=
-(-1)^{p_ip_{i+1}}x_1\wedge\ldots\wedge x_{i+1}\wedge x_{i}\wedge\ldots
\wedge x_{n-1}$.)
Conversely, if $(L,V,ad)$ is a triple with $L$ a Lie superalgebra, $V$ an
$L$-module, and $ad$ an $L$-module morphism from $\wedge^{n-1}V \mapsto L$
 such that the map
 $(v_1,\ldots,v_n)\mapsto  ad(v_1\wedge\ldots\wedge v_{n-1})\, v_n$
 from $\otimes^{n}V$ to $V$ is $\mathbb{Z}_2$-skewsymmetric
(we call the homomorphisms of this type {\it skewsymmetric}),
 then $V$ becomes an  $n$-ary Filippov superalgebra by defining
\begin{center} $[v_1,\ldots,v_n]=ad(v_1\wedge\ldots\wedge v_{n-1})\, v_n.$ \end{center}
Therefore, we obtain a correspondence between the set of $n$-ary
 Filippov superalgebras and the set of triples $(L,V,ad)$,
 satisfying the conditions above.

We assume that all vector spaces appearing in the following are
finite-dimensional over $F$.

If $\F$ is a simple $n$-ary Filippov superalgebra then Theorem \ref{teo1}
shows that the Lie superalgebra $Inder(\F)$  is semisimple, and  $\F$
 is a faithful and irreducible $Inder(\F)$-module. Moreover,
 the $Inder(\F)$-module morphism  $ad:\wedge^{n-1}\F\mapsto Inder(\F)$
 is surjective.

Conversely, if $(L,V,ad)$ is a triple such that $L$ is a semisimple
 Lie superalgebra over $F$, $V$ is a
faithful irreducible $L$-module, $ad$ is a surjective
 $L$-module morphism from $\wedge^{n-1}V$ onto the adjoint module $L$,
 and the map
 $(v_1,\ldots,v_n)\mapsto ad(v_1\wedge\ldots\wedge v_{n-1})\, v_n$
 from $\otimes^{n}V$ to $V$ is $\mathbb{Z}_2$-skewsymmetric, then the corresponding
 $n$-ary Filippov superalgebra is simple. A triple with these conditions
 will be called a {\it good triple}. Thus, the problem of determining
 the simple  $n$-ary Filippov superalgebras over $F$ can be
 translated to that of finding the good triples.

\section{Some notations and results on Lie superalgebras}

In this section, we recall some notations and results from \cite{Kac}
on the Lie superalgebra $A(m,n)$ (and its irreducible faithful
finite-dimensional representations). We also give some explicit constructions
which we shall use some later in the study of the simple finite-dimensional Filippov superalgebras of type $A(m,n)$.
Let us start recalling the definition
of induced module.

Let $\L$ be a Lie superalgebra, $U(\L)$ its universal enveloping superalgebra
\cite{Kac}, $H$ a subalgebra of $\L$, and $V$ an $H$-module. The module
$V$ can be extended to $U(H)$-module. We consider the $\mathbb{Z}_2$-graded space
$U(\L)\otimes_{U(H)}V$ (this is the quotient space of $U(\L)\otimes V$ by the
linear span of the elements of the form $gh\otimes v-g\otimes h(v)$,
$g\in U(\L)$, $h\in U(H)$). This space can be endowed with the structure of a
$\L$-module as follows $g(u\otimes v)=gu\otimes v,g\in \L,u\in U(\L),v\in V$.
The so-constructed $\L$-module is said to be {\it induced from the $H$-module
$V$} and is denoted by $Ind_H^{\L}V$.

From now on, we denote by $G$  a contragredient Lie superalgebra over
$\Phi$ and  consider it with the ``standard'' $\mathbb{Z}$-grading
\cite[Sections 5.2.3 and 2.5.7]{Kac}.

Let $G=\oplus_{i\geq -d}G_i$.
 Set $H=(G_0)\0=\left< h_1,\ldots,h_n\right>$, $N^+=\oplus_{i>0}G_i$ and
 $B=H\oplus N^+$. Let $\Lambda\in H^*, \Lambda(h_i)=a_i\in\Phi$,
$\left<v_\Lambda \right>$ be an one-dimensional $B$-module for which
$N^+(v_\Lambda)=0, h_i(v_\Lambda)=a_iv_\Lambda$. Let $\delta_i \in H^*$, $\delta_i(h_j)=\delta_{ij}$ where
$\delta_{ij}$ is Kronecker's delta.
Let $V_\Lambda=Ind_B^G\left<v_\Lambda \right>/I_{\Lambda}$,
where $I_{\Lambda}$ is the (unique) maximal submodule of the $G$-module
$V_\Lambda$. Then $\Lambda$ is called
the {\it highest weight} of the $G$-module
 $V_\Lambda$.
 By \cite{Kac}, every faithful irreducible finite-dimensional
 $G$-module may be obtained in this manner. Note that the condition
$1\otimes v_\Lambda\in V\0(V\1)$ gives a $\mathbb{Z}_2$-grading on $V_\Lambda$.

\begin{lemma} \textnormal{\cite{Poj8}} \label{soma} \it
 Let $V$ be a module over a Lie superalgebra $G$,
 let $V=\oplus V_{\gamma_i}$ be its weight decomposition, and let $\phi$ be
 a homomorphism from $\wedge^mV$ into $G$. Then,  for all $v_i\in
 V_{\gamma_i}$,

\begin{center}
\begin{tabular} {ll}
 $\phi(v_1,\ldots,v_m)\in G_{\gamma_1+\ldots+\gamma_m}$,
& \hspace{0,5cm} if \ $\gamma_1+\ldots+\gamma_m$ is a root of $G$, \\
$\phi(v_1,\ldots,v_m)=0$, & \hspace{0,5cm} otherwise.
\end{tabular}
\end{center}
 \end{lemma}

Let $G$ be a contragredient Lie superalgebra of rank $n$,
$U=Ind_B^G\left<v_\Lambda \right>$, and
$V=V_\Lambda=U/N$
be a finite-dimensional representation of $G$,
where $N=I_{\Lambda}$ is a maximal proper submodule of the $G$-module
$V_\Lambda$.
Let $G=\oplus_{\alpha}G_{\alpha}$
 be a root decomposition of $G$ relative to a Cartan subalgebra $H$.
Denote by ${\cal A}$ the following set of roots$:$
${\cal A}=\{\alpha:g_{\alpha}\notin B\}$.

\begin{lemma} \textnormal{\cite{PS}} \label{pesos} \it Let $g_{\alpha}\in G_{\alpha}$ and $g_{\alpha}\otimes
v\neq 0$ $(v=v_\Lambda)$. Then

\begin{center} $g_{\alpha}^j\otimes v\in U_{\sum_{i=1}^n
(j\alpha(h_i)+\Lambda(h_i))\delta_i}$ \end{center}
for all $j\in{\mathbb N}$, and there
exists a minimal
 $ k\in {\mathbb N}$ such that  $g_{\alpha}^k\otimes v \in N$. Moreover, the set
${\cal E}_{\alpha,k}=\{1\otimes v,g_{\alpha}\otimes v,\ldots,
g_{\alpha}^{k-1}\otimes v\}$ is linearly independent in $V$.
Setting $h=[g_{-\alpha},g_{\alpha}]$, we have

\begin{enumerate}
\item $\Lambda(h)=-\frac{(k-1)\alpha(h)}{2}$ if either $g_{\alpha}\in G\0$
   or $k\notin 2{\mathbb N}$;

\item $ \alpha(h)=0$ if $g_{\alpha}\in G\1$ and $k\in 2{\mathbb N}$.
\end{enumerate}
 \end{lemma}

\begin{remark}
Note that if we start with a root $\beta$ then there exists
$s\in{\mathbb N}$ such that ${\cal E}_{\beta,s}$ is linearly independent,
but ${\cal E}_{\alpha,k}\cup {\cal E}_{\beta,s}$
may not be linearly independent.
\end{remark}

Recall that a set ${\cal E}$ is called a {\it pre-basis} of a vector space $W$
if $\left<{\cal E}\right>=W$.

Let $\{g_{\alpha_1}^{k_1}\ldots g_{\alpha_s}^{k_s}\otimes v; k_i\in {\mathbb
N}_0, \alpha_i\in{\cal A}\}$ be a pre-basis of $U$. As we have seen above, for
every $i=1,\ldots,s$, there exists a minimal $p_i\in {\mathbb N}$ such that
$g_{\alpha_i}^{p_i}\otimes v \in N$. Using the induction on the word length, it is easy
to show that $\{g_{\alpha_1}^{k_1}\ldots g_{\alpha_s}^{k_s}\otimes v; k_i\in
{\mathbb N}_0,k_i<p_i,\alpha_i\in{\cal A}\}$ is a pre-basis of $U/N$.

We finish this part with some more notations that we use in the two next sections:

\hspace{-0,7cm} $\bullet$ the symbol $\= $
denotes an equality up to a nonzero coefficient;

\hspace{-0,7cm} $\bullet$ $\underline{u,v}_{\ t}$ means that the elements $u$ and
$v$ are $t$-times repeating $\underbrace{u,v,\ldots,u,v}_{2t}$\,, being
the index $t$ omitted when its value is clear from the context.

\section{Simple Filippov superalgebras of type $A(m,n)$}

In what follows, considering $A(m,n)$, we assume that $m \neq n$. Recall that $A(m,n):=sl(m+1,n+1)$ for $m \neq n$ and $m, n \in \mathbb{N}_0$. It consists of the matrices of type
$$\left(
\begin{array}{c|c}
A         &   B \\
\hline
&\\[-3mm]
C    &  D\\
\end{array}
\right),
$$
 where $A \in M_{(m+1)\times(m+1)}(F), B \in M_{(m+1)\times (n+1)}(F), C \in M_{(n+1) \times (m+1)}(F), D \in M
 _{(n+1)\times (n+1)}(F)$ and $tr(A)=tr(D)$. Let us write some elements in $G=A(m,n)$:
\begin{eqnarray*}
&&\left.
\begin{array}{l}
 h_i =\ e_{ii}-e_{i+1,i+1}, \ \ i=1,\ldots,m,m+2,\ldots,m+n+1,\\
 h_{m+1}=e_{m+1,m+1}+e_{m+2,m+2},\\
 e_{kl}:=g_{\epsilon_k-\epsilon_l} \in G_{\epsilon_k-\epsilon_l}, \ \ k,l=1,\ldots,m+1 \ \textnormal{or} \ k,l=m+2,\ldots,m+n+2,
\end{array}
\right\}\in G\0\\
&&\left.
\begin{array}{l}
e_{kl}:=g_{\epsilon_k-\epsilon_l} \in G_{\epsilon_k-\epsilon_l}, \ \ k=1,\ldots,m+1, l=m+2,\ldots,m+n+2, \\
\hspace{3,6cm} \textnormal{or} \
k=m+2,\ldots,m+n+2, l=1,\ldots,m+1.
\end{array}
\right\}\in G\1
\end{eqnarray*}

The space $H:=G_0=\left< h_1,\ldots,h_{m+n+1}\right>$
is a Cartan subalgebra
of $A(m,n)$, and $\epsilon_i$
are the linear functions
on $H$ defined by its values on $h_1,\ldots,h_{m+n+1}$
and the conditions $\epsilon_i(e_{jj})=\delta_{ij}$, where $\delta_{ij}$
is Kronecker's delta. Then $\Delta=\Delta_0\cup\Delta_1$ is
a root system for $A(m,n)$, where $\Delta_0=\{ 0; \epsilon_k-\epsilon_l, k,l=1,\ldots,m+1 \ \textnormal{or} \ k,l=
m+2,\ldots,m+n+2\}$,
and $\Delta_1=\{\epsilon_k-\epsilon_l, k=1,\ldots,m+1, l=m+2, \ldots,m+n+2 \ \textnormal{or} \ k=m+2,\ldots,m+n+2, l=1,\ldots,m+1\}$.
The roots $\{\alpha_i:= \epsilon_i-\epsilon_{i+1},i=1,\ldots,m+n+1\}$
 are simple.

The conditions $deg \ g_{\alpha_i}=1, deg \ g_{-\alpha_i}=-1$ give us the standard grading of $A(m,n)$, \cite[Section 5.2.3]{Kac}. The negative part of this grading is $G_{\epsilon_k-\epsilon_l}, l < k$. Because of this, the set
\begin{eqnarray}\label{exp}{\cal E}=\left\{
\prod_{l < k} g_{\epsilon_k-\epsilon_l}^{\gamma_{kl}}
\otimes v\ :\ \gamma_{kl} \in{\mathbb N}_0\right\}
\end{eqnarray}
is a pre-basis of the induced module
$M=Ind_B^G\left<v_\Lambda \right>$ (where $v=v_\Lambda$).

Let $V$ be an irreducible module over $G=A(m,n)$ with the highest
weight $\Lambda$, $\Lambda(h_i)=a_i$. Denote $\Lambda$ by $(a_1,\ldots,a_{m+n+1})$.
Applying Lemma \ref{pesos}, we have

\begin{equation} a_i \in \mathbb{N}_0 \ \ \textnormal{if} \ i \neq m+1. \label{am+1}\end{equation}

If $u\in V_{\gamma}$ (or $G_{\gamma}$) then we may write $h_q u=p_q (u)u$
($[h_q ,u]=p_q (u)u$) and call {\it $q$-weight} of $u$ to $p_q(u)$.

In what follows, the symbol $w \oti{i}{j}$ means that
$p_{h_l}(w \otimes v)=p_l(w \otimes v)=i$ and $p_{h_r}(w \otimes v)=p_r(w \otimes v)=j$ (the same with the notation
$\mathop u\limits_{j}^{i}$), where everytime in text we specify the indices $l$ and $r$.

\subsection{The type $A(1,n)$}

In this subsection, because of \cite{BPoj} and \cite{pat_artigo4}, we assume that $n \in \mathbb{N} \setminus \{1\}$.
We begin with some technical
lemmas on irreducible modules of some special types over $A(1,n)$.

\begin{lemma} \label{miercoles}
Let $V=V_{\Lambda}$ be an irreducible module over $G=A(1,n)$ with $\Lambda=(a_1, \ldots, a_{n+2})$,
$a_1=1$ and $a_{2} \neq 0$. Assume that $(G,V,\phi)$ is a good triple. Then $a_{2} \leq -\frac 1 2$.
\end{lemma}

\proof Suppose that $\phi(u_1,\ldots,u_s)=g_{\epsilon_{3}-\epsilon_2}$. Consider $h=h_1+h_{2}=e_{11}+e_{33}$.
From the nonzero action on $1 \otimes v$, we obtain $|2-p_{h}(u_i)+a_{2}| \leq 1$. If $a_{2} > -\frac 1 2$
then $p_{h}(u_i) > \frac 1 2$, which leads to a contradiction since $p_{h}(g_{\epsilon_{3}-\epsilon_2})=1$. Thus,
$a_{2} \leq -\frac 1 2$. \qed

\begin{lemma} \label{piatnitsa2}
Let $V=V_{\Lambda}$ be an irreducible module over $A(1,n)$ with $\Lambda=(a_1, \ldots, a_{n+2})$,
$a_1=1$ and $a_2=\ldots=a_{2+n}=0$. Then

\begin{equation*} \left\{\prod_{i=2}^{3+n} g_{\epsilon_i-\epsilon_1}^{\alpha_i} \otimes v: \alpha_i
\in \{0,1\}\right\} \end{equation*} is a pre-basis of $V$.
\end{lemma}

\proof Consider $g_{\epsilon_i-\epsilon_1}$ with $i \neq 1$. Suppose that $g_{\epsilon_i-\epsilon_1} \in G_{\0}$.
Then, by Lemma \ref{pesos}, $g_{\epsilon_i-\epsilon_1}^2 \otimes v = 0$. If $g_{\epsilon_i-\epsilon_1} \in G_{\1}$
then, from $[g_{\epsilon_i-\epsilon_1}, g_{\epsilon_i-\epsilon_1}]\otimes v=0$, we obtain the same conclusion. \qed

\begin{lemma} Let $V=V_{\Lambda}$ be an irreducible module over $A(1,n)$ with
$\Lambda=(a_1, \ldots, a_{n+2})$, $a_2=a\neq 0$ and $a_i=0$ for $i \neq 2$.
Suppose that $h=e_{11}+e_{33}$ and $\phi(u_1,\ldots,u_s)=g_{\epsilon_3-\epsilon_2}$. Then it is impossible that $p_h(u_i)=a$ for all $i \in \{1,\ldots,s\}$.
\label{lemaquase}
\end{lemma}

\proof Notice that $w_0=\phi(1 \otimes v , u_2, \ldots, u_s) \in {\cal E}_1$,
where ${\cal E}_1$ denotes the elements from $A(1,n)$ with $h$-weight equal to $1$:
 ${\cal E}_1=\{g_{\epsilon_1-\epsilon_i}, g_{\epsilon_3-\epsilon_i}: i \neq 1, 3\}$.
If $w_0= g_{\epsilon_3-\epsilon_i}$ then we can multiply $w_0$ by $g_{\epsilon_i-\epsilon_{n+3}}$
to think that $w_0=g_{\epsilon_3-\epsilon_{n+3}}$. We may proceed analogously with $g_{\epsilon_1-\epsilon_i}$.
 Now, if $w_0=g_{\epsilon_1-\epsilon_{n+3}}$ then we may multiply it by
$g_{\epsilon_3-\epsilon_1}$ to arrive at $g_{\epsilon_3-\epsilon_{n+3}}$.
Thus, we may replace all $u_i$ either with $g_{\epsilon_{j_i}-\epsilon_2} \otimes v\ (j_i>3)$ or $ 1 \otimes v$
(maybe, multiplied by $\alpha:=g_{\epsilon_3-\epsilon_1}$).  Thus, we arrive at
$$w=\phi(\alpha^{\delta_1}g_{\epsilon_{i_1}-\epsilon_2} \otimes v, \ldots,
\alpha^{\delta_r}g_{\epsilon_{i_r}-\epsilon_2} \otimes v,
\alpha^{\delta_{r+1}}g_{\epsilon_{n+3}-\epsilon_2} \otimes v, \ldots,\qquad \qquad \qquad \qquad$$
$$\qquad \qquad \qquad \qquad \qquad \qquad \qquad \qquad \ldots,
\alpha^{\delta_{r+q}}g_{\epsilon_{n+3}-\epsilon_2} \otimes v,\underline{\alpha\otimes v}_{t},
\underline{1\otimes v}_{p}) \=g_{\epsilon_3-\epsilon_{n+3}},$$
where $3<i_j<n+3, \delta_k\in \{0,1\}.$
The action of $h$ on $w$ gives $a(r+q+t+p)=1$, and the action of $e_{22}+e_{n+3,n+3}$ gives $r=2$.
The action of $e_{n+2,n+2}-e_{n+3,n+3}$ gives either $q=1$ or $q=0$.

If $\mathbf{ q=1}$ then $i_1=i_2=n+2$, and the action of $e_{33}-e_{n+2,n+2}$ implies $t+\sum_{i=1}^3\delta_i=3$.
If $\delta_3=0$ then $g_{\epsilon_{n+3}-\epsilon_2} \otimes v\curvearrowright g_{\epsilon_{n+3}-\epsilon_1} \otimes v$
 gives an element with $(e_{11}+e_{n+3,n+3})$-weight being equal to $-2$.
 If $\delta_1=0$ (or $\delta_2=0$) then
 $g_{\epsilon_{n+2}-\epsilon_2} \otimes v\curvearrowright g_{\epsilon_{n+3}-\epsilon_1} \otimes v$
leads to an element with $(e_{11}+e_{n+2,n+2})$-weight being equal to $-2$.
 Thus, $\delta_1=\delta_2=\delta_3=1, t=0$. Consider the action of
 $w$ on $g_{\epsilon_{n+3}-\epsilon_1} \otimes v$ with the consecutive change
 $\alpha g_{\epsilon_{n+3}-\epsilon_2} \otimes v\curvearrowright g_{\epsilon_{n+3}-\epsilon_1} \otimes v$.
 We get $$w_1=\phi(\alpha g_{\epsilon_{n+2}-\epsilon_2} \otimes v,\alpha g_{\epsilon_{n+2}-\epsilon_2} \otimes v,
  g_{\epsilon_{n+3}-\epsilon_1} \otimes v,\underline{1\otimes v}_p),$$
and $p_1(w_1)=-1,p_2(w_1)=p_{n+2}(w_1)=1, $ i.e., $w_1\=g_{\epsilon_2-\epsilon_{n+3}}$.
The action of $w_1$ on $g_{\epsilon_{n+3}-\epsilon_2} \otimes v$ gives
$$w_2=\phi(\alpha g_{\epsilon_{n+2}-\epsilon_2} \otimes v,g_{\epsilon_{n+3}-\epsilon_2} \otimes v,
  g_{\epsilon_{n+3}-\epsilon_1} \otimes v,\underline{1\otimes v}_p)\neq 0$$
  (note that we may assume $n=2$, since otherwise the $(e_{22}+e_{44})$-weight of $w$ is $-2$).
  Then $p_1(w_2)=p_2(w_2)=p_3(w_2)=0, p_4(w_2)=-1$, and there is no an element
   in $A(1,2)$ with such root.

  If $\mathbf{q=0}$ then $i_1=n+2, 3<i_2<n+3$, and the action of $e_{33}-e_{n+2,n+2}$ implies $t+\sum_{i=1}^2\delta_i=2$.
If $\delta_1=0$ then $g_{\epsilon_{n+2}-\epsilon_2} \otimes v\curvearrowright g_{\epsilon_{n+3}-\epsilon_1} \otimes v$
 gives a $(e_{11}+e_{n+2,n+2})$-contradiction. Moreover, $p_{e_{22}+e_{44}}(w)=-1$ if $i_2\neq 4$.
 If $n>3$ then the action of $e_{22}+e_{55}$ gives a contradiction.
 Thus, we may assume $i_2=4,n=3$ and
$$w=\phi(\alpha g_{\epsilon_{5}-\epsilon_2} \otimes v,\alpha^{\delta_2} g_{\epsilon_{4}-\epsilon_1} \otimes v,
\underline{1\otimes v})\= g_{\epsilon_{3}-\epsilon_6}.$$
The action on $g_{\epsilon_{6}-\epsilon_1}$
 gives
$$w_1=\phi(\alpha g_{\epsilon_{5}-\epsilon_2} \otimes v,g_{\epsilon_{6}-\epsilon_1} \otimes v,
\underline{1\otimes v})\neq 0.$$ We have
$p_1(w_1)=-1=p_4(w_1),p_2(w_1)=p_{3}(w_1)=1,p_{5}(w_1)=0 $, i. e., $w_1\=g_{\epsilon_2-\epsilon_{4}}$.
From the action on $ g_{\epsilon_{4}-\epsilon_2} \otimes v$, we arrive at
$$w_2=\phi(g_{\epsilon_{4}-\epsilon_2} \otimes v,g_{\epsilon_{6}-\epsilon_1} \otimes v,
\underline{1\otimes v})\neq 0.$$
We have $p_1(w_2)=p_2(w_2)=0,p_3(w_2)=-1=p_5(w_2), p_4(w_2)=1$, and there is no an element
   in $A(1,3)$ with such root.
\qed

\begin{lemma} Let $V=V_{\Lambda}$ be an irreducible module over $A(1,n)$ with
$\Lambda=(a_1, \ldots, a_{n+2})$, $a_2=a\neq 0$ and $a_i=0$ for $i \neq 2$. Suppose that $h=e_{11}+e_{33}$. Then it is impossible to have

\begin{equation} \phi(\stackrel{1+a}{u_1}, \ldots, \stackrel{1+a}{u_k}, \stackrel{a}{v_1}, \ldots, \stackrel{a}{v_{s-k}})=g_{\epsilon_3-\epsilon_2}, \label{ig}\end{equation} where the superscripts denote the $h$-weights, for some $k \geq 2$.
\label{lemafinal}
\end{lemma}

\proof Denote by $\mathcal{E}_r$ the set of elements from $A(1,n)$ with $h$-weight equal to $r$ .
Suppose first that $s-k>0$. From (\ref{ig}), acting on $1 \otimes v$, we have $w_1=\phi(u_1,\ldots,u_k, 1
\otimes v, v_2, \ldots, v_{s-k}) \in \mathcal{E}_{1}=\{g_{\epsilon_1-\epsilon_i}, g_{\epsilon_3-\epsilon_i}:
i \neq 1, 3\}$. Assume that $1 \otimes v$ is odd. We may think that $w_1 \neq g_{\epsilon_3-\epsilon_2}$ since,
 otherwise, we can multiply $w_1$ by $g_{\epsilon_2-\epsilon_i} (i >3)$. If $w_1=g_{\epsilon_3-\epsilon_i} (i >3)$
 then, through the multiplication by $g_{\epsilon_1-\epsilon_3}$, we arrive at $w_1'=g_{\epsilon_1-\epsilon_i}$. Thus,
 $w_1 \= g_{\epsilon_1-\epsilon_i}$ and we may act on $t_1=g_{\epsilon_j-\epsilon_1} \otimes v$ ($j=i$ if $i>3$;
 $i\neq j>3$ if $i=2$), interchanging this element with $u_1$. We obtain
 $w_2=\phi(\stackrel{a-1}{t_1},\stackrel{a+1}{u_2}, \ldots, \stackrel{a+1}{u_k},\stackrel{a}{v_1}, \ldots,
 \stackrel{a}{v_{s-k}}) \in \mathcal{E}_{-1}=\{g_{\epsilon_i-\epsilon_1}, g_{\epsilon_i-\epsilon_3}: i \neq 1,3\}$,
where, here and throughout this proof, the $h$-weights are above the elements.
 If $w_2=g_{\epsilon_i-\epsilon_1}$ then we may multiply it by $g_{\epsilon_1-\epsilon_3}$ to arrive at
 $g_{\epsilon_i-\epsilon_3}$. So, $w_2 \= g_{\epsilon_i -\epsilon_3}$ and we replace $u_2$ by the action on
 $t_2=g_{\epsilon_3-\epsilon_2} \oti{a+1}{}$. Repeating this procedure, we substitute $t_1$ and arrive at a
 skewsymmetry contradiction. Suppose that $1 \otimes v$ is even.
We have $w_1= g_{\epsilon_3-\epsilon_i}$ $(i >3)$ or $w_1=g_{\epsilon_1-\epsilon_j}$ $(j \neq 1, 3)$.
In the former case we may multiply $w_1$ by $g_{\epsilon_1-\epsilon_3}$ to get the latter one. In the latter case we
can act on $t_1$ to arrive at $w_2$. Here we repeat the above argument to arrive at
$w_3=\phi(\underline{t_2}, v_1, \ldots, v_{s-k}) \in \mathcal{E}_1$. If $w_3=g_{\epsilon_3-\epsilon_i} (i > 3)$
then we may multiply $w_3$ by $g_{\epsilon_i-\epsilon_{n+3}}$ to think that $w_3=g_{\epsilon_3-\epsilon_{n+3}}$.
We can do the same with $g_{\epsilon_{1}-\epsilon_i} (i>3)$. Thus, we may assume that
$\mathcal{E}_1=\{g_{\epsilon_3-\epsilon_2},g_{\epsilon_3-\epsilon_{n+3}},g_{\epsilon_1-\epsilon_2},
g_{\epsilon_1-\epsilon_{n+3}}\}$. Notice that if $w_3=g_{\epsilon_1-\epsilon_{n+3}}$ then we can multiply it by
$g_{\epsilon_3-\epsilon_1}$ to arrive at $g_{\epsilon_3-\epsilon_{n+3}}$. Let $h_0=e_{22}+e_{n+3,n+3}$. We remark
that the element $g_{\epsilon_3-\epsilon_1}$ does not change either the $h$-weights or the $h_0$-weights. Therefore,
we may replace all $v_i$ with $g_{\epsilon_i-\epsilon_2} \otimes v, g_{\epsilon_{n+3}-\epsilon_2} \otimes v, 1 \otimes v$
(maybe multiplied by $g_{\epsilon_3-\epsilon_1}$). Adding the $h_0$-weights, we get $-2k-r=-2$ for some
 $r \in \mathbb{N}_0$, which is impossible because $k \geq 2$.

Now suppose that $s-k=0$. Thus, (\ref{ig}) has the following shape

\begin{equation} \phi(\stackrel{1+a}{u_1}, \ldots, \stackrel{1+a}{u_s})=g_{\epsilon_3-\epsilon_2}. \end{equation} We may multiply it by $g_{\epsilon_2-\epsilon_i}$ $(i > 3)$ to assume that $\phi(u_1,\ldots,u_s)= g_{\epsilon_3-\epsilon_i}$. Through the multiplication by $g_{\epsilon_1-\epsilon_3}$, we may assume that $\phi(u_1,\ldots,u_s)=g_{\epsilon_1-\epsilon_i}$. We can now interchange $u_1$ and $g_{\epsilon_i-\epsilon_1} \oti{a-1}{}$, and repeat the above described procedure to substitute all $u_i$ with $t_2=g_{\epsilon_3-\epsilon_2} \otimes v$. Thus, we arrive at $\phi(\underline{t_2}) \in \mathcal{E}_1$. Considering the $3$-weights, we obtain $\phi(t_2,t_2)=g_{\epsilon_3-\epsilon_4}$. From here, thinking in the $1$-weights, we have a weight contradiction.
\qed

\begin{lemma} \label{lim}
Let $V=V_{\Lambda}$ be an irreducible module over $G=A(1,n)$ with $\Lambda=(a_1, \ldots, a_{n+2})$.
Suppose that $a_1=0, a_2 \neq 0$ and $\sum_{i=3}^{n+2} a_i= 1$. Assume that $(G,V,\phi)$ is a good triple.
 Then $0 < a_2 \leq 1/2$.
\end{lemma}

\proof Suppose that $\phi(u_1,\ldots,u_s)=g_{\epsilon_3-\epsilon_2}$ and consider $H=h_1+h_2=e_{11}+e_{33}$.
By the nonzero action over $1 \otimes v$, we obtain $|1-p_H(u_i)+a_2| \leq 1$. From here and taking
into account that $\sum_{i=1}^s p_H(u_i)=1$, we conclude that $a_2 \leq 1/2$. Let $H'=e_{22}+e_{n+3,n+3}$.
Assume that $\phi(u_1,\ldots,u_s)=g_{\epsilon_{n+3}-\epsilon_2}$. By the action on $1 \otimes v$,
we have $|-p_{H'}(u_i)+a_2-1| \leq 1$. Whence, $a_2 >0$. \qed

\begin{theorem} There are no simple finite-dimensional Filippov
superalgebras of type $A(1,n)$ over $\Phi$. \end{theorem}

\proof Suppose that $V$ is a finite-dimensional irreducible module
over $G=A(1,n)$ with the highest weight $\Lambda=(a_1, \ldots, a_{2+n})$ ($a_1\neq 0$),
and
 $\phi$ is a surjective skewsymmetric homomorphism from $\wedge^sV$ on $G$. Then there exist $u_i\in V_{\gamma_i}$ such that
\begin{eqnarray}
  \label{eqnova}
  \phi(u_1,\ldots,u_{s})&=&g_{\epsilon_{2}-\epsilon_1}.
\end{eqnarray} By Lemma \ref{soma}, $\sum_{i=1}^sp_1(u_i)=-2$. From Lemma \ref{pesos},
  $g_{\epsilon_{2}-\epsilon_1}^{a_1}\otimes v\neq 0$.
 Since $\phi$
is a skewsymmetric homomorphism,
$\phi(u_1,\ldots,u_{i-1},g_{\epsilon_{2}-\epsilon_1}^{a_1-1}\otimes v,
u_{i+1},\ldots,u_{s})\neq 0$.
As $p_1(g_{\epsilon_{2}-\epsilon_1}^{a_1-1}\otimes v)=2-a_1$, the inequality
$|p_1(u_i)+a_1|\leq 2$ follows.  From here we see that
the required skewsymmetric homomorphism does not exist if $a_1\geq 4$.

From now on, unless stated otherwise, we put the $1$-weights above the elements.

Consider the case \framebox{$a_1=3$}. In this case, $p_1(u_i) <0$.
So, by (\ref{eqnova}), we have $\phi(\stackrel{-1}{u_1},\stackrel{-1}{u_2})=g_{\epsilon_{2}-\epsilon_1}$ and, acting on $1 \otimes v$, we arrive at $\phi(\stackrel{-1}{u_1}, 1\oti{3}{})\=
g_{\epsilon_1-\epsilon_{2}}$. Acting on $g_{\epsilon_{2}-\epsilon_1}\otimes v$, we
obtain $\phi(g_{\epsilon_{2}-\epsilon_1}\oti{1}{ },1\oti{3}{ })\neq 0$,
which is a weight contradiction.

Now let us take \framebox{$a_1=2$}. As $\sum_{i=1}^s p_1(u_i)=-2$ and, in this case, $p_1(u_i) \leq 0$, we can only have

\begin{center} i) $\phi(\stackrel{-2}{u_1}, \stackrel{0}{u_{2}}, \ldots,  \stackrel{0}{u_s})
=g_{\epsilon_{2}-\epsilon_1}$ \ \ \ or \ \ ii) $\phi(\stackrel{-1}{u_1}, \stackrel{-1}{u_2},
\stackrel{0}{u_3}, \ldots, \stackrel{0}{u_s})
=g_{\epsilon_{2}-\epsilon_1}$. \end{center}
First consider i). Let us suppose that $1 \otimes v$ is even.
Acting on $1 \otimes v$,
we have $\phi(1 \oti{2}{},
\ox{u_2}{0}{},\ldots,\ox{u_{s}}{0}{})\=
g_{\epsilon_1-\epsilon_{2}}$.
Then, acting twice on $g_{\epsilon_{2}-\epsilon_1} \otimes v$, we arrive at
$\phi(\underline{g_{\epsilon_{2}-\epsilon_1} \otimes v}_{2}, u_3, \ldots, u_s) \neq 0$ which leads to a skewsymmetry
contradiction.
To finish the consideration of this subcase, suppose now that $1 \otimes v$  is odd. Then, acting on $1 \otimes v$
and, repeatedly, on $g_{\epsilon_{2}-\epsilon_1}\otimes v$, we get
$\phi(\stackrel{2}{1 \otimes v},\underline{g_{\epsilon_{2}-\epsilon_1} \stackrel{0}{\otimes} v})
\=g_{\epsilon_1-\epsilon_{2}}$. From here, analizing the $2$-weights, we
conclude that $a_{2}=-1$. Assume that $\phi(u_1,\ldots,u_s)=g_{\epsilon_{3}-\epsilon_{2}}$.
Consider $h'=h_1+h_2=e_{11}+e_{33}$. From the nonzero action
on $1 \otimes v$, we have $|-p_{h'}(u_i)+2| \leq 1$. So, we obtain a contradiction because $\sum_{i=1}^s p_{h'}(u_i)=1$.
In the case ii), the multiplication
by $g_{\epsilon_1-\epsilon_{2}}$
gives either $$\phi(\stackrel{-1}{u_1}, \stackrel{-1}{u_2},
\stackrel{0}{u_3}, \ldots, \stackrel{2}{w},\ldots, \stackrel{0}{u_s})
(1\otimes v)\neq 0\  \mbox{ or }\
\phi(\stackrel{-1}{u_1}, \stackrel{1}{v_2},
\stackrel{0}{u_3}, \ldots, \stackrel{0}{u_s})(1\otimes v)\neq 0,$$
for some $w, v_2$.
In both cases, replacing $u_1$ by $1\otimes v$, we arrive at a
weight contradiction.

Now take \framebox{$a_1=1$}. Consider (\ref{eqnova}) and
$h''=h_{3}+\ldots+h_{2+n}$. By the nonzero action on $1 \otimes v$, we have $|-p_{h''}(u_i)+a_{3}+\ldots+a_{2+n}| \leq 2$. So, we can deduce that $a_{3}+\ldots+a_{2+n} <3$ since otherwise it is impossible to have (\ref{eqnova}). As $a_r \in \mathbb{N}_0$ for $r \neq 2$ then

\begin{center} $a_{3}+\ldots+a_{2+n} \in \{0,1,2\}$.\end{center} In what follows,
we analise these three possibilities, numbered with I), II) and III),
for $a_{3}+\ldots+a_{2+n}$.

I) Assume that $a_{3}+\ldots+a_{2+n}=2$ and let $h=h_1+h_{3}+\ldots+h_{2+n}$.
 From (\ref{eqnova}), by the action on $1 \otimes v$, we arrive at $|1-p_{h}(u_i)| \leq 2$.
 Thus, $p_{h}(u_i) \in \{-1,0,1,2,3\}$ and we obtain $\phi(\stackrel{-1}{u_1}, \stackrel{-1}{u_2},
\stackrel{0}{\overline{u_3, \ldots, u_s}})
= g_{\epsilon_{2}-\epsilon_1}$, where the $h$-weights are above the elements.
The multiplication by $g_{\epsilon_{1}-\epsilon_{2}}$ leads to

\begin{center}
$\phi(\stackrel{-1}{u_1}, \stackrel{1}{w},
\stackrel{0}{\overline{u_3, \ldots, u_s}})(1 \oti{3}{}) \neq 0$ or $\phi(\stackrel{-1}{u_1}, \stackrel{-1}{u_2},
 \stackrel{2}{\overline{\ldots}})
(1 \oti{3}{}) \neq 0$,
\end{center} for some $ w$. In both cases, replacing $u_1$ by $1 \otimes v$, we get a weight contradiction.

II) Take $h$ as above and suppose that $a_{3}+\ldots+a_{2+n}=1$.
Once again by the action of $g_{\epsilon_{2}-\epsilon_1}$ on $1 \otimes v$, we get $|p_{h}(u_i)| \leq 2$.
So, we have either $\phi(\stackrel{-1}{u_1}, \stackrel{-1}{u_2},
\stackrel{0}{\overline{u_3, \ldots, u_s}}) = g_{\epsilon_{2}-\epsilon_1}$ or $
\phi(\stackrel{-2}{u_1},
\stackrel{0}{\overline{u_2, \ldots, u_s}}) =g_{\epsilon_{2}-\epsilon_1}$. In the former subcase,
using the reasoning of I), we obtain more weight contradictions. In the latter subcase, multiplying by
$g_{\epsilon_1-\epsilon_2}$, we arrive either at $\phi(\stackrel{-2}{u_1},
 \stackrel{2}{\overline{\ldots}})(1 \oti{2}{}) \neq 0$, which gives a
weight contradiction, or at
$\phi(\stackrel{0}{w}, \stackrel{0}{\overline{u_2, \ldots, u_s}})(1 \oti{2}{}) \neq 0$. Thus, we have $\phi(1 \oti{2}{}, \stackrel{0}{\overline{u_2, \ldots, u_s}}) \in \{g_{\epsilon_1-\epsilon_2}, g_{\epsilon_1-\epsilon_{n+3}}, g_{\epsilon_3-\epsilon_2}, g_{\epsilon_3-\epsilon_{n+3}}\}$. Taking into account that $p_{h}(g_{\epsilon_3-\epsilon_1})=0=p_{h}(g_{\epsilon_{n+3}-\epsilon_2})$ and making adequate multiplications, we may assume that $\phi(\stackrel{2}{z}, \stackrel{0}{\overline{u_2, \ldots, u_s}}) \= g_{\epsilon_3-\epsilon_2}$, for some $z$. Suppose first that $g_{\epsilon_3-\epsilon_2} \otimes v \neq 0$. Then, by the action on $1 \oti{2}{}$, we get a weight contradiction. Suppose now that $g_{\epsilon_3-\epsilon_2} \otimes v = 0$.
Thus $a_2=0$.  Assume that $\phi(u_1, \ldots, u_s)=g_{\epsilon_{n+3}-\epsilon_2}$. Let $h=e_{22}+e_{n+3,n+3}$ and, in what follows, consider the $2$-weights above the elements and the $h$-weights underneath them. From the action on $1 \otimes v$, it is possible to conclude that
$p_2(u_i), p_h(u_i) \in \{-2,-1,0\}$. Thus we have $\phi(\mathop {u_1}\limits_{0}^{-1},\mathop {u_2}\limits_{0}^{0},\ldots,\mathop {u_s}\limits_{0}^{0}) = g_{\epsilon_{n+3}-\epsilon_2}$. From here, considering the $2$-weights and the $h$-weights, we arrive at

\begin{center} $\phi(\mathop {u_1}\limits_{0}^{-1},1\oti{0}{-1},\mathop {u_3}\limits_{0}^{0},\ldots,\mathop {u_s}\limits_{0}^{0}) \= g_{\epsilon_i-\epsilon_2}$, where $i \neq 2, 3, n+3$. \end{center} If
$g_{\epsilon_i-\epsilon_2} \otimes v \neq 0$ then the action on $1 \oti{}{-1}$ gives a weight contradiction. So,
$g_{\epsilon_i-\epsilon_2} \otimes v = 0$ and there exists $j>i$ such that $g_{\epsilon_i-\epsilon_2}g_{\epsilon_j-\epsilon_i} \otimes v \neq 0$. If $j \neq n+3$ then $p_h(g_{\epsilon_j-\epsilon_i} \otimes v)=-1$ and we obtain a weight contradiction through the action on $g_{\epsilon_j-\epsilon_i} \oti{}{-1}$. Thus we may assume that $a_{n+2}=1$. Let us replace all $u_k \ (k \geq 3)$ with
$g_{\epsilon_{n+3}-\epsilon_i} \otimes v$ and act one more time on such element. We have

\begin{center} $\phi(1 \oti{0}{-1}, g_{\epsilon_{n+3}-\epsilon_{i_1}}\oti{0}{0}, \ldots,
g_{\epsilon_{n+3}-\epsilon_{i_{s-1}}}\oti{0}{0}) \neq 0$. \end{center} Considering the
$h$, $2$ and $1$-weights, we conclude that

\begin{equation}
\phi(1 \oti{0}{-1}, g_{\epsilon_{n+3}-\epsilon_1}\oti{0}{0}, \ldots,
g_{\epsilon_{n+3}-\epsilon_{1}}\oti{0}{0}) \= g_{\epsilon_1-\epsilon_{n+3}}.
\label{Santander1}
\end{equation} From the multiplication by
$g_{\epsilon_{n+3}-\epsilon_1}$, we have
$\phi(\underline{g_{\epsilon_{n+3}-\epsilon_1}} \otimes v) \= h_1+h_2-h_3- \ldots - h_{n+2}$.
The action on $1 \otimes v$ leads to a contradiction.

III)
We now have $a_t=0$ for $t \in \{3, \ldots, 2+n\}$. Suppose that $a_{2}>0$. Consider $H'= h_1 + h_{2}= e_{11}+e_{33}$ and assume that

\begin{equation} \phi(u_1,\ldots,u_s)=g_{\epsilon_{3}-\epsilon_1}. \label{Akademgorodok1} \end{equation}
Acting on $1 \otimes v$, we have $|-p_{H'}(u_i)+1+a_{2}|\leq 1$ and, consequently, $p_{H'}(u_i) \geq a_{2} > 0$,
being (\ref{Akademgorodok1}) impossible. Thus, $a_{2} \leq 0$.

Let us take first $a_{2} < 0$. Suppose that

\begin{equation} \phi(u_1,\ldots,u_s)=g_{\epsilon_{3}-\epsilon_{2}}. \label{Akademgorodok2} \end{equation}
By the action on $1 \otimes v$, we arrive at $p_{2}(u_i) \leq 1 +a_{2}$. We can't have (\ref{Akademgorodok2})
if $a_{2} < -1$. So, we conclude that $a_{2} \geq -1$. Now assume that
$\phi(u_1,\ldots,u_s)=g_{\epsilon_{3+n}-\epsilon_{2}}$. Consider the action on $1 \otimes v$. On one hand,
we have $|-1-p_{2}(u_i)+a_{2}| \leq 1$ and, so, $a_2-2 \leq p_{2}(u_i) \leq a_2 < 0$. On the other hand,
$|2-p_1(u_i)| \leq 2$. Thus, we have $p_1(u_i) \in \{0,1\}$. From now on, in this subcase and unless stated otherwise,
we will put the $1$-weights above the elements and the $2$-weights underneath them. Taking into account the $1$-weights,
we can only have

\begin{equation} \phi(\stackrel{1}{u_1},\stackrel{0}{u_2},\ldots,\stackrel{0}{u_s})=g_{\epsilon_{3+n}-\epsilon_{2}}.
\label{Akademgorodok3}  \end{equation} Acting on $1 \otimes v$ allows us to obtain
$\phi(\stackrel{1}{u_1},1 \oti{1}{a_{2}},\stackrel{0}{u_3},\ldots,\stackrel{0}{u_s})\=g_{\epsilon_{1}-\epsilon_{2}}$.
Notice that $a_{2}$ has to be greater than $-1$; otherwise, we obtain a weight contradiction.
 By Lemma \ref{miercoles}, $a_{2} \leq -\frac 1 2$. Therefore, we arrive at $a_2=-1/2$ and
 $\phi(\mathop {u_1}\limits_{-1/2}^{1},1 \oti{1}{-1/2})\=g_{\epsilon_{1}-\epsilon_{2}}$. From
 the action on $g_{\epsilon_{3+n}-\epsilon_1} \otimes v$, we have

\begin{equation} \phi(g_{\epsilon_{3+n}-\epsilon_1}\oti{0}{-1/2},1 \oti{1}{-1/2}) \neq 0. \label{Akademgorodok4}
\end{equation} Taking into account the $1$-weights, the $2$-weights and the $(2+n)$-weights in (\ref{Akademgorodok4}),
we conclude that we must have

\begin{equation} \phi(g_{\epsilon_{n+3}-\epsilon_1} \otimes v,1 \otimes v) \= g_{\epsilon_{n+3}-\epsilon_{2}}.
\label{Akademgorodok5}  \end{equation} The nonzero action on $1 \otimes v$ leads to

\begin{equation} \phi(1 \oti{1}{-1/2},1 \oti{1}{-1/2}) \= g_{\epsilon_{1}-\epsilon_{2}}. \label{X}\end{equation}
If $1 \otimes v$ is even then we have a skewsymmetry contradiction. If $1 \otimes v$ is odd then we obtain
the contradiction $0=g_{\epsilon_1-\epsilon_{3}}$ from the multiplication by $g_{\epsilon_{2}-\epsilon_{3}}$
in (\ref{X}).

Now let $a_{2}=0$. Observe that, by Lemma \ref{piatnitsa2},
all $2$-weights of the elements of the pre-basis of $V$ are zero or positive. Therefore, it is impossible to find $v_i \in V$ such that $\phi(v_1,\ldots,v_s)=g_{\epsilon_1-\epsilon_{2}}$.

Now take \framebox{$a_1=0$}. Suppose first that $a_2 \neq 0$.
Assume that $\phi(u_1,\ldots,u_s)=g_{\epsilon_3-\epsilon_1}$ and $H''=h_3+\ldots+h_{n+2}=e_{33}-e_{n+3,n+3}$.
By the action on $1 \otimes v$, we have $|1-p_{H''}(u_i)+a_3+\ldots+a_{n+2}| \leq 2$. From here, as
$\sum_{i=1}^s p_{H''}(u_i)=1$, we conclude that $a_3+\ldots+a_{n+2} <2$. Taking into account (\ref{am+1}),
 we have to study two subcases: 1) $a_3=\ldots=a_{n+2}=0$;
2) $\sum_{i=3}^{n+2}a_i=1$.

1) Consider $a_2=a$, $h=e_{11}+e_{33}, v_i=g_{\epsilon_i-\epsilon_1}, w_j=g_{\epsilon_j-\epsilon_2}$ and

\begin{equation} \phi(u_1,\ldots,u_s)=g_{\epsilon_3-\epsilon_2}. \label{quase} \end{equation}
Then $\mathcal{E} =\left< v_{i_1}\ldots v_{i_t}w_{j_1}\ldots w_{j_r} \otimes v: i_p, j_q \neq 1, 2 \right>$
is a pre-basis of $V$. Note that $p_h(u_i)=a + t_i$, where $t_i \in \mathbb{Z}$ and $t_i \leq 1$. Then
$\sum_{i=1}^s (t_i+a)=1$ and, by the action of (\ref{quase}) on $1 \otimes v$, we have $|1-(t_j+a)+a| \leq 1$
for all $j$. Thus, $t_j \in \{0, 1\}$ and $k(1+a)+(s-k)a=1$ for some $k \in \mathbb{N}_0$. By Lemma \ref{lemaquase},
$k \neq 0$. Notice also that $k \neq 1$.
So, $k \geq 2$. We can also see that $-1 < a <0$. From Lemma \ref{lemafinal}, we have that, for $k \geq 2$, this subcase can not occur.

2) Consider $\phi(u_1, \ldots, u_s)=g_{\epsilon_3-\epsilon_2}$ and $H=e_{22}+e_{n+3,n+3}$. From the action
on $1 \otimes v$, we obtain $|-1-p_H(u_i)+a_2-1| \leq 1$. From here and by Lemma \ref{lim}, we get
$-3 < p_H(u_i) \leq -\frac 1 2$. Therefore, $a_2=\frac 1 2, s=2$ and $p_H(u_1)=p_H(u_2)=-\frac 1 2$.
Notice that, from the same action, for $h=e_{33}-e_{n+3,n+3}$, we get $p_h(u_i) \geq 0$. Henceforth, we have
 \begin{center} $\phi(\mathop {u_1}\limits_{1}^{-1/2},\mathop {u_2}\limits_{0}^{-1/2})=g_{\epsilon_3-\epsilon_2}$,
 \end{center} where
the $H$-weights are above the elements and the $h$-weights underneath them.
Acting on $1 \otimes v$, we obtain
$\phi(\mathop {u_1}\limits_{1}^{-1/2},1 \oti{-1/2}{1}) \=g_{\epsilon_{3}-\epsilon_{n+3}}$. The action on
 $g_{\epsilon_{n+3}-\epsilon_2} \otimes v$, taking into account the weights over $H, h, e_{11}+e_{n+3,n+3}$ and $h_3$,
leads to $\phi(g_{\epsilon_{n+3}-\epsilon_2} \otimes v, 1\otimes v) \= g_{\epsilon_1-\epsilon_{n+3}}$.
Through the multiplication by
$g_{\epsilon_{n+3}-\epsilon_2}$ we deduce that $1\otimes v$ is even and we arrive at $\phi(
\underline{g_{\epsilon_{n+3}-\epsilon_2} \otimes v}) \= g_{\epsilon_1-
\epsilon_2}$. From here, multiplying by $g_{\epsilon_2-\epsilon_1}$, we have the contradiction
$0 = h_1$.

 At last, suppose that $a_2 = 0$. We may assume that $\sum_{i=3}^{n+2} a_i=a >0$. Consider
 $\phi(u_1,\ldots,u_s)=g_{\epsilon_{n+3}-\epsilon_2}$. The action on $1 \otimes v$ leads
 to $|-1-p_2(u_i)| \leq 1$ and, consequently, $p_2(u_i) \leq 0$ for all $i$. If $h=e_{22}+e_{n+3,n+3}$
 then $|-p_h(u_i)-a| \leq 1$. So, $a=1$ and $p_h(u_i)=0$ for all $i$. Therefore, we have
 $\phi(\mathop {u_1}\limits_{0}^{-1},1\oti{0}{-1},\mathop {u_3}\limits_{0}^{0},\ldots,\mathop {u_s}\limits_{0}^{0})
  = g_{\epsilon_{i}-\epsilon_2}$ where $i \neq 2, 3, n+3$, the $2$-weights and the $h$-weights
  are above and underneath the elements, respectively. If $g_{\epsilon_i-\epsilon_2} \otimes v \neq 0$
  then $u_1 \mapsto 1 \otimes v$ gives a $h$-weight contradiction. If $g_{\epsilon_i-\epsilon_2} \otimes v = 0$
  then there exists a $j>i$ such that $g_{\epsilon_i-\epsilon_2}g_{\epsilon_j-\epsilon_i} \otimes v \neq 0$.
  If $j \neq n+3$ then $p_h(g_{\epsilon_j-\epsilon_i} \otimes v)=-1$ and $u_1 \mapsto g_{\epsilon_j-\epsilon_i} \otimes v$
  gives a $h$-weight contradiction. Thus, we may assume that $a_{n+2}=1$. We can replace $u_k$ with
  $g_{\epsilon_{n+3}-\epsilon_i} \otimes v$, $k\geq 3$. Continuing the process, we obtain
  $w_1=\phi(1\otimes v, g_{\epsilon_{n+3}-\epsilon_{i_1}} \otimes v,
  \ldots, g_{\epsilon_{n+3}-\epsilon_{i_{s-1}}} \otimes v) \neq 0$. Since the $h_{n+2}$-weight of
   $w_1$ is $1$, $w_1 \in \{g_{\epsilon_j-\epsilon_{n+3}}, g_{\epsilon_{n+2}-\epsilon_j}: j \neq n+2, n+3\}$.
   Let $h_0=e_{11}+e_{n+3,n+3}$. But $p_{h_0}(w_1)=-1$. Hence, either $w_1 \= g_{\epsilon_{n+2}-\epsilon_1}$,
   which gives a $h_1$-weight contradiction, or $w_1\=g_{\epsilon_i-\epsilon_{n+3}}$. In the latter case,
   considering the $h_k$-weights for $k=1,2,3,\ldots$, we arrive at a weight contradiction.

To finish the proof, consider now $a_i=0$ for all $i$. Then $V$ is trivial. \qed

\subsection{The main theorem}

We can now state and prove the main result of this article.

\begin{theorem} There exist no simple finite-dimensional Filippov
superalgebras of type $A(m,n)$ over $\Phi$. \label{mainresult} \end{theorem}

\proof We can suppose that $G=A(m,n)$ with $m \neq n$ and $m \geq 2$, because we have already proved
  that there exist no simple Filippov superalgebras of type $A(n,n)$ with $n \in \mathbb{N}$, \cite{BPoj}, nor of type $A(1,n)$ with
$n \in \mathbb{N}_0\setminus \{1\}$ and of type $A(0,n)$ with $n \in \mathbb{N}$, \cite{pat_artigo4}.
Assume that $V$ is a finite-dimensional irreducible module
over $G$ with the highest weight $\Lambda=(a_1, \ldots, a_{m+n+1})$ ($a_1+\ldots+a_m\neq 0$),
and
 $\phi$ is a surjective skewsymmetric homomorphism from $\wedge^sV$ on $G$.
Then there exist $u_i\in V_{\gamma_i}$ such that
\begin{eqnarray}
  \label{eq:l261}
  \phi(u_1,\ldots,u_{s})&=&g_{\epsilon_{m+1}-\epsilon_1}.
\end{eqnarray} Let $H=h_1+\ldots+h_m=e_{11}-e_{m+1,m+1}$. By Lemma \ref{soma}, $\sum_{i=1}^sp_H(u_i)=-2$. From Lemma \ref{pesos},
  $g_{\epsilon_{m+1}-\epsilon_1}^{a_1+\ldots+a_m}\otimes v\neq 0$.
 Since $\phi$
is a skewsymmetric homomorphism,
$\phi(u_1,\ldots,u_{i-1},g_{\epsilon_{m+1}-\epsilon_1}^{a_1+\ldots+a_m-1}\otimes v,
u_{i+1},\ldots,u_{s})\neq 0$.
As $p_H(g_{\epsilon_{m+1}-\epsilon_1}^{a_1+\ldots+a_m-1}\otimes v)=2-a_1-\ldots-a_m$, the inequality
$|p_H(u_i)+a_1+\ldots+a_m|\leq 2$ follows.  From here we see that
the required skewsymmetric homomorphism does not exist if $a_1+\ldots+a_m\geq 4$.

Throughout this proof, unless stated otherwise, we put the $H$-weights above the elements.

Consider the case \framebox{$a_1+\ldots+a_m=3$}. Then we have $\phi(\stackrel{-1}{u_1},\stackrel{-1}{u_2})=g_{\epsilon_{m+1}-\epsilon_1}$ and, ac\-ting on $1 \otimes v$, we arrive at $\phi(\stackrel{-1}{u_1}, 1\oti{3}{})\=
g_{\epsilon_1-\epsilon_{m+1}}$. By the action on $g_{\epsilon_{m+1}-\epsilon_1}\otimes v$, we
obtain $\phi(g_{\epsilon_{m+1}-\epsilon_1}\oti{1}{ },1\oti{3}{ })\neq 0$,
which is a weight contradiction.

Now let us take \framebox{$a_1+\ldots+a_m=2$}. Thus, there are two possibilities

\begin{center} i) $\phi(\stackrel{-2}{u_1}, \stackrel{0}{u_{2}}, \ldots,  \stackrel{0}{u_s})
=g_{\epsilon_{m+1}-\epsilon_1}$ \ \ \ or \ \ ii) $\phi(\stackrel{-1}{u_1}, \stackrel{-1}{u_2},
\stackrel{0}{u_3}, \ldots, \stackrel{0}{u_s})
=g_{\epsilon_{m+1}-\epsilon_1}$. \end{center}
First consider i). Let us suppose that $1 \otimes v$ is even.
Acting on $1 \otimes v$,
we have $\phi(1 \oti{2}{},
\ox{u_2}{0}{},\ldots,\ox{u_{s}}{0}{})\=
g_{\epsilon_1-\epsilon_{m+1}}$.
Then, acting twice on $g_{\epsilon_{m+1}-\epsilon_1} \otimes v$, we arrive at
$\phi(\underline{g_{\epsilon_{m+1}-\epsilon_1} \otimes v}_{2}, u_3, \ldots, u_s) \neq 0$ which leads to a skewsymmetry
contradiction.
To fi\-nish the consideration of this subcase,
suppose now that $1 \otimes v$  is odd. Then, acting on $1 \otimes v$ and, repeatedly, on
$g_{\epsilon_{m+1}-\epsilon_1}\otimes v$, we get
$\phi(\stackrel{2}{1 \otimes v},\underline{g_{\epsilon_{m+1}-\epsilon_1}\stackrel{0}{\otimes}v})
\=g_{\epsilon_1-\epsilon_{m+1}}$. From here, analyzing the $(m+1)$-weights, we
conclude that $a_{m+1}=-1$. Assume that \begin{equation} \phi(u_1,\ldots,u_s)=g_{\epsilon_{m+2}-\epsilon_{m+1}}. \label{mais}\end{equation} From the nonzero action
on $1 \otimes v$, we arrive at $|3-p_{H}(u_i)| \leq 2$. Consequently, we can't have (\ref{mais}). In the case ii), the multiplication
by $g_{\epsilon_1-\epsilon_{m+1}}$
gives, for some $v_i$, either $$\phi(\stackrel{-1}{u_1}, \stackrel{-1}{u_2},
\stackrel{0}{u_3}, \ldots, \stackrel{2}{v_i},\ldots, \stackrel{0}{u_s})
(1\oti{2}{})\neq 0\  \mbox{ or }\
\phi(\stackrel{-1}{u_1}, \stackrel{1}{v_2},
\stackrel{0}{u_3}, \ldots, \stackrel{0}{u_s})(1\oti{2}{})\neq 0, \ \textnormal{for some $v_i, v_2$}.$$
In both cases, replacing $u_1$ by $1\otimes v$, we arrive at a
weight contradiction.

Now consider \framebox{$a_1+\ldots+a_m=1$}. Suppose that $a_{m+1} >0$. Take
 $h=h_1+...+h_{m+1}=e_{11}+e_{m+2, m+2}$, and assume that

\begin{equation} \phi(u_1,\ldots,u_s)=g_{\epsilon_{m+2}-\epsilon_1}. \label{Santander3}\end{equation} Through the nonzero action on $1 \otimes v$, we have $p_h(u_i) > 0$ and (\ref{Santander3}) can not occur. Thus, $a_{m+1} \leq 0$.

I) Suppose that $a_{m+1} <0$. Consider $\phi(u_1,\ldots,u_s)=g_{\epsilon_{m+2}-\epsilon_{m+1}}$.
By the action on $1 \otimes v$,
taking into account the $(m+1)$-weights and the $h$-weights, we arrive at $-1 \leq a_{m+1} \leq -\frac 1 2$.

Ia) Assume that $a_{m+1} \neq -1$. Let $\phi(u_1,\ldots,u_s)=g_{\epsilon_{m+2}-\epsilon_{m+1}}$. By
the action on $1 \otimes v$, we have $p_{H}(u_i) \in \{0,1,2,3,4\}$. Consequently, after the
action on $1 \otimes v$, we arrive at

\begin{center} $\phi(\stackrel{1}{u_1}, 1 \oti{1}{}, \stackrel{0}{u_3},\ldots,\stackrel{0}{u_s})
\=g_{\epsilon_{1}-\epsilon_{m+1}}$.\end{center} Replacing every $u_k (k \geq 3)$ by $g_{\epsilon_{m+2}-\epsilon_1} \otimes v$ and acting one more time on the mentioned element, we get
$\phi(1 \oti{1}{}, g_{\epsilon_{m+2}-\epsilon_1} \oti{0}{}, \ldots, g_{\epsilon_{m+2}-\epsilon_1}
\oti{0}{}) \neq 0$.
Analyzing the $H$, $h$ and $1$-weights involved, we conclude that $a_1=1$ and

\begin{center} $\phi(1 \oti{1}{}, g_{\epsilon_{m+2}-\epsilon_1} \oti{0}{}, \ldots, g_{\epsilon_{m+2}-\epsilon_1} \oti{0}{}) \= g_{\epsilon_1-\epsilon_i}, \ i\neq 1, 2, m+1, m+2$. \end{center}
Through the multiplication by $g_{\epsilon_{m+2}-\epsilon_1}$, we obtain $\phi(\underline{g_{\epsilon_{m+2}-\epsilon_1} \otimes v}) \= g_{\epsilon_{m+2}-\epsilon_i}, i \neq 1, 2, m+1, m+2$.
Thus, considering the $2$-weights, we have $i \neq 3$. Continuing the process, through
the consecutive analises of the $2, 3, \ldots$-weights, we eliminate all the possibilities for $i$.

Ib) Assume that $a_{m+1}=-1$. $1$) Consider $a_1=1$ and suppose that
$\phi(u_1,\ldots, u_s)=g_{\epsilon_{m+2}-\epsilon_{m+1}}$. In this subcase,
we put the $(m+1)$-weights above the elements and the $1$-weights underneath them.
Through the action on $1 \otimes v$, we have $p_{m+1}(u_i) \in \{-2,-1,0\}, p_H(u_i) \in \{0,1,2,3,4\}$ and $p_1(u_i) \in \{-1,0,1,2,3\}$. If there is a $k$ such that $p_1(u_k)=-1$ then the replacement of $u_k$ by $1 \otimes v$ leads to a $(m+1), 1$-weights contradiction.
Thus, $p_1(u_i) \geq 0$ and, through the action
of $g_{\epsilon_{m+2}-\epsilon_{m+1}}$ on $1 \otimes v$, putting the $H$-weights in the third line, we arrive at

\begin{equation*} \phi(\mathop{\mathop{u_1}\limits_{0}^{0}}\limits_1,
\mathop{1 \oti{-1}{1}}\limits_{1},
\mathop{\mathop{u_3}\limits_0^0}\limits_{0},\ldots,\mathop{\mathop{u_s}\limits_0^0}\limits_{0})\= g_{\epsilon_{1}-\epsilon_{m+1}}. \end{equation*} By the action on $g_{\epsilon_{m+1}-\epsilon_1} \otimes v$, we have
$\phi(\mathop{g_{\epsilon_{m+1}-\epsilon_1}\oti{0}{0}}\limits_{\hspace{0,9cm} -1},
\mathop{1 \oti{-1}{1}}\limits_{1},
\mathop{\mathop{u_3}\limits_0^0}\limits_{0},\ldots,\mathop{\mathop{u_s}\limits_0^0}\limits_{0}) \neq 0$.
This is a weight contradiction since we don't have an element in $A(m,n)$ with the obtained $m+1, 1, H$-weights.
$2$) Now consider $a_m=1$ and suppose that $\phi(u_1,\ldots,u_s)=g_{\epsilon_{m+2}-\epsilon_{m+1}}$. By the action
on $1 \otimes v$, we conclude that $p_h(u_i) \in \{0,1,2\}$, $p_{m+1}(u_i) \in \{-2,-1,0\}$ and
$p_m(u_i) \in \{0,1,2,3,4\}$. So, we have

\begin{equation} \phi(\mathop {u_1}\limits_{0}^{1},\mathop {u_2}\limits_{0}^{0},\ldots,
\mathop {u_s}\limits_{0}^{0})= g_{\epsilon_{m+2}-\epsilon_{m+1}}, \label{Labra1} \end{equation}
where the $h$-weights are above the elements and the $(m+1)$-weights are underneath them.
From the action on $1 \otimes v$, we obtain

\begin{center} $\phi(1 \mathop {\otimes}\limits_{-1}^{0} v,\mathop {u_2}\limits_{0}^{0},\ldots,\mathop
{u_s}\limits_{0}^{0}) \= g_{\epsilon_1-\epsilon_{m+2}} (g_{\epsilon_{i}-\epsilon_{m+1}}, i \neq 1, m+1, m+2)$.
\end{center} Consider the former possibility. As $p_m(g_{\epsilon_{m+2}-\epsilon_{m+1}})=1$,
$p_m(g_{\epsilon_1-\epsilon_{m+2}})=0$, $p_m(1 \otimes v)=1$ then, in (\ref{Labra1}), $p_m(u_1)=2$ and
the sum of the $m$-weights of the remaining elements is equal to $-1$, which is impossible. Consider the latter occasion.
Notice that $p_m(g_{\epsilon_i-\epsilon_{m+1}})$ is either $1$ (when $i \neq m$) or $2$ (when $i=m$). If $i \neq m$
then $p_m(u_1)=1$ and the change $u_2 \mapsto 1 \otimes v$ in (\ref{Labra1}) gives a $h, m+1, m$-weights contradiction.
 If $i=m$ then $p_m(u_1)=0$ and, acting on $1 \otimes v$ in (\ref{Labra1}), we have
 $\phi(\mathop{1 \oti{0}{-1}}\limits_1,\mathop{\underline{\mathop {u_2}\limits_{0}^{0},\ldots,
\mathop {u_s}\limits_{0}^{0}}}\limits_1) \= g_{\epsilon_m-\epsilon_{m+1}}$, where the $m$-weights of
the elements are in the third weight line. Through the action on
$g_{\epsilon_{m+1}-\epsilon_m} \otimes v$, we obtain $\phi(1 \mathop{\oti{0}{-1}}\limits_1, \mathop{g_{\epsilon_{m+1}-\epsilon_m}\oti{0}{0}}\limits_{\hspace{1cm}-1},
\mathop{\mathop{u_3}\limits_0^0}\limits_0,\ldots,\mathop{\mathop{u_s}\limits_0^0}\limits_0) \neq 0$, one more weight
contradiction. $3$) Assume that there exists a $j \in \{2, \ldots, m-1\}$ such that $a_j=1$. In
this subcase, we put the $(m+1)$-weights above the elements and the $j$-weights underneath them.
Assume that $\phi(u_1, \ldots, u_s)=g_{\epsilon_{j+1}-\epsilon_j}$. The action on $1 \otimes v$
allows us to conclude that $p_{m+1}(u_i) \in \{-2,-1,0\}$ and $p_j(u_i) \in \{-3,-2,-1,0,1\}$. If, for example,  $p_j(u_1) = 1$ then, through the mentioned action, we obtain the weight contradiction $\phi(1 \mathop {\otimes}\limits_{1}^{-1} v,\mathop\protect{\underline{\stackrel{0}{u_2},\ldots,\stackrel{0}{u_s}}}\limits_{-3}) \neq 0$. Suppose now that $\phi(\mathop {u_1}\limits_{0}^{0},\ldots,\mathop {u_{s-2}}\limits_{0}^{0},\mathop {u_{s-1}}\limits_{-1}^{0},\mathop {u_{s}}\limits_{-1}^{0}) = g_{\epsilon_{j+1}-\epsilon_{j}}$. Multiplying the last equality by $g_{\epsilon_j-\epsilon_{j+1}}$ and acting on $1 \otimes v$, we get $\phi(\mathop\protect{\underline{\stackrel{0}{u_1},\ldots,\stackrel{0}{u_{s-1}}}}\limits_{1},1 \oti{-1}{1}) \neq 0$, which is a weight contradiction. Finally, we study the subcase $u=\phi(\mathop {u_{1}}\limits_{0}^{0},\ldots,\mathop {u_{s-1}}\limits_{0}^{0},\mathop {u_{s}}\limits_{-2}^{0})=g_{\epsilon_{j+1}-\epsilon_j} (*)$. The action on $1 \otimes v$ leads to $w=\phi(\mathop {u_{1}}\limits_{0}^{0},\ldots,\mathop {u_{s-1}}\limits_{0}^{0},1 \oti{-1}{1}) \in \{g_{\epsilon_j-\epsilon_{m+1}},g_{\epsilon_j-\epsilon_{m+2}}\}$. As $p_{j-1}(u)=1$ and $p_{j-1}(w)=-1$ then $p_{j-1}(u_s)=2$. So, through the action on $1 \otimes v$ in $(*)$, we obtain the weight contradiction
$\phi(\mathop{1 \oti{-1}{1}}\limits_{0},
\mathop{\underline{\mathop {u_2}\limits_{0}^{0},\ldots,\mathop {u_{s-1}}\limits_{0}^{0}}}\limits_{\geq 0},
\mathop{\mathop{u_s}\limits_{-2}^0}\limits_{2}) \neq 0$, where the third weight line refers to $h_{j-1}$.

II) Consider the case $a_{m+1}=0$. Assume that
$\phi(u_1,\ldots, u_s)=g_{\epsilon_{m+2}-\epsilon_1}.$ Take $h$ as above.
Through the action on $1 \otimes v$, we have $p_{m+1}(u_i), p_h(u_i) \in \{0,1,2\}$. Thus

\begin{equation} \phi(\mathop {u_1}\limits_{0}^{1},\mathop {u_2}\limits_{0}^{0},
\ldots,\mathop {u_s}\limits_{0}^{0})= g_{\epsilon_{m+2}-\epsilon_1}, \label{Cantabria}
\end{equation} where, here and in what follows, we consider the $(m+1)$-weights above the
elements and the $h$-weights underneath them. By the action on $1 \otimes v$ in (\ref{Cantabria}), we arrive at

\begin{equation} \phi(\mathop {u_1}\limits_{0}^{1},1 \oti{0}{1},\mathop {u_3}\limits_{0}^{0},\ldots,
\mathop {u_s}\limits_{0}^{0})=g_{\epsilon_{m+2}-\epsilon_i}, \textnormal{with $ i \neq 1, m+1, m+2$}.
\label{Cantabria1} \end{equation}
If $g_{\epsilon_{m+2}-\epsilon_i} \otimes v \neq 0$ then we get a weight contradiction from the action on
$1 \oti{}{1}$. So, $g_{\epsilon_{m+2}-\epsilon_i} \otimes v =0$ and there exists a $j < i$ such that
$g_{\epsilon_{m+2}-\epsilon_i}g_{\epsilon_i-\epsilon_j} \otimes v \neq 0$. If $j \neq 1$ then
$p_h(g_{\epsilon_i-\epsilon_j} \otimes v)=1$ and from the action on $g_{\epsilon_i-\epsilon_j} \otimes v$
arises a weight contradiction. We may assume that $a_1=1$. Let us replace all $u_k$ in (\ref{Cantabria1}),
for $k \geq 3$,
by $g_{\epsilon_{i_k}-\epsilon_1} \oti{0}{0}$ and act one more time on
$g_{\epsilon_{i_1}-\epsilon_1} \oti{0}{0}$. Then, looking at the $(m+1), h,
1$-weights, we obtain
\begin{center}
$u:=\phi(1 \otimes v, g_{\epsilon_{i_1}-\epsilon_1} \otimes
v,g_{\epsilon_{i_3}-\epsilon_1} \otimes v,
\ldots, g_{\epsilon_{i_{s}}-\epsilon_1} \otimes v)=g_{\epsilon_1-\epsilon_t}
(g_{\epsilon_{m+2}-\epsilon_{m+1}})$,
\end{center} with $t \neq 1,2,m+1, m+2$.
Note that $i_k<m+1$, since otherwise $q:=i_k=m+2$ for some $k$ and we arrive at
$(e_{11}+e_{qq})$-contradiction.
Moreover, if $u=g_{\epsilon_{m+2}-\epsilon_{m+1}}$ then the multiplication on
$g_{\epsilon_1-\epsilon_{m+2}}$ gives a contradiction. Now, for $t =m, m-1, \ldots, 2$,
considering, consecutively, all these $t$-weights, we arrive at
a weight contradiction.

Finally, suppose that \framebox{$a_1+\ldots+a_m=0$}.
As $A(m,n) \simeq A(n,m)$, \cite[Section 4.2.2]{Kac}, then $a_{m+2}+\ldots+a_{m+n+1}=0$. Thus, consider
$a_t=0$ for $t\neq m+1$ and $a_{m+1}=a \neq 0$. Assume that $h=e_{11}+e_{m+n+2,m+n+2}$. Let $w=\phi(u_1,\ldots, u_s)=g_{\epsilon_{m+1}-\epsilon_2}$. Then $p_h(w)=0$ and $p_h(u_i)=a+\epsilon_i$ with $\epsilon_i \in \mathbb{Z}$.
Take $x=g_{\epsilon_{m+n+2}-\epsilon_{m+1}} \otimes v$.
We have $wx \neq 0$ and $p_h(x)=a+1$. If $p_h(u_i)=a-\epsilon_i$, for some $i$ and $\epsilon_i \in \mathbb{N}$,
then $u_i \curvearrowright x$ gives a weight contradiction. Therefore, $a <0$. Now let
$w=\phi(u_1,\ldots,u_s)=g_{\epsilon_{m+3}-\epsilon_{m+2}}$. Then $p_h(w)=0$ and $p_h(u_i)=a + \epsilon_i$
with $\epsilon_i \in \mathbb{Z}$. Take $x=g_{\epsilon_{m+2}-\epsilon_1} \otimes v$. Then $p_h(x)=a-1$ and $wx \neq 0$.
If $p_h(u_j)=a+\epsilon_j$, for some $j$ and $\epsilon_j \in \mathbb{N}$, then $u_j \curvearrowright x$ gives a weight
contradiction. Henceforth, $a >0$. Thus, $a=0$ and the module is trivial. This
finishes the proof of the theorem. \qed

\begin{corollary}
There is no simple finite-dimensional Filippov superalgebra ${\cal F}$ of type
$A(m,n)$ such that ${\cal F}$ is a highest weight module over $A(m,n)$.
\end{corollary}

\renewcommand{\refname}{\begin{center} References \end{center}}

\vspace{0,8cm}

\begin{footnotesize}P. D. Beites \\ Departamento de Matem\'{a}tica and Centro de Matemсtica,
Universidade da Beira Interior \\ Covilhу, Portugal
 \\ \textit{E-mail adress}: pbeites@ubi.pt\end{footnotesize}

\vspace{0,4cm}

\begin{footnotesize}A. P. Pozhidaev \\ Sobolev Institute of Mathematics and Novosibirsk State University \\
Novosibirsk, Russia \\ \textit{E-mail adress}: app@math.nsc.ru \end{footnotesize}
\vspace{.5cm}

\end{document}